\documentclass[12pt]{amsart}

\date{\today}

\usepackage{amsmath,amsthm,amssymb,amsfonts,amscd, array}
\usepackage[margin=1in]{geometry}
\usepackage{setspace}
\usepackage[hang,flushmargin,symbol*]{footmisc}
\usepackage{mathtools}
\usepackage{enumitem}
\usepackage{hyphenat}
\usepackage{graphicx}
\usepackage{subfig}
\usepackage{color}
\definecolor{darkblue}{rgb}{0, 0, .6}
\definecolor{grey}{rgb}{.7, .7, .7}
\usepackage{hyperref}
\hypersetup{
	colorlinks=true,
	linkcolor=darkblue,
	anchorcolor=darkblue,
	citecolor=darkblue,
	pagecolor=darkblue,
	urlcolor=darkblue,
	pdftitle={},
	pdfauthor={}
}

\theoremstyle{definition}
\newtheorem{thm}{Theorem}[subsection] 

\newtheorem{lem}[thm]{Lemma}
\newtheorem{prop}[thm]{Proposition}
\newtheorem{defn}[thm]{Definition}
\newtheorem{rem}[thm]{Remark}
\newtheorem{ex}[thm]{Example}

\newcommand{\N}{\mathbb{N}}

\newcommand{\C}{\widetilde{C}}
\renewcommand{\O}{\mathcal{O}}
\newcommand{\E}{\mathcal{E}}
\newcommand{\z}{\mathsf{z}}
\newcommand{\x}{\mathsf{x}}
\newcommand{\TL}{TL}
\newcommand{\supp}{\mathrm{supp}}
\renewcommand{\L}{\mathcal{L}}
\newcommand{\R}{\mathcal{R}}
\renewcommand{\(}{\left(}
\renewcommand{\)}{\right)}
\newcommand{\w}{\mathsf{w}}
\renewcommand{\H}{\mathcal{H}}
\newcommand{\FC}{\mathrm{FC}}
\renewcommand{\r}{\mathbf{r}}

\title[Non-cancellable elements in type affine $C$ Coxeter groups]{Non-cancellable elements in type \\ affine $C$ Coxeter groups}

\author{Dana C.~Ernst}
\address{Department of Mathematics and Statistics, Northern Arizona University, Flagstaff, AZ} 
\email{\url{dana.ernst@nau.edu}}
\urladdr{\url{http://danaernst.com}}


\begin{document}

\begin{abstract}
Let $(W,S)$ be a Coxeter system and suppose that $w \in W$ is fully commutative (in the sense of Stembridge) and has a reduced expression beginning (respectively, ending) with $s \in S$.  If there exists $t\in S$ such that $s$ and $t$ do not commute and $tw$ (respectively, $wt$) is no longer fully commutative, we say that $w$ is left (respectively, right) weak star reducible by $s$ with respect to $t$.  We say that a fully commutative element is non-cancellable if it is irreducible under weak star reductions.  In this paper, we classify the non-cancellable elements in Coxeter groups of types $B$ and affine $C$.  In a sequel to this paper, the classification of the non-cancellable elements play a pivotal role in inductive arguments used to prove the faithfulness of a diagrammatic representation of a generalized Temperley--Lieb algebra of type affine $C$.
\end{abstract}

\maketitle

%%%%%%%%%%%%%%%%%%%%%

\begin{section}{Introduction}\label{sec:intro}

Let $(W,S)$ be a Coxeter system with group $W$ and finite set of generating involutions $S$.  Our principal focus in this paper will be the infinite Coxeter group of type affine $C$, denoted $\C$, and its finite subgroups of type $B$.  A well-known result in the theory of Coxeter groups, known as Matsumoto's Theorem, states that any two reduced expressions for $w \in W$ are equivalent under the equivalence relation generated by braid relations. If $w$ is such that any two of its reduced expressions are equivalent by iterated commutations of commuting generators, $w$ is called fully commutative~\cite{Stembridge.J:B}.  We denote the set of fully commutative elements of $W$ by $\FC(W)$.  Fully commutative elements arise in several contexts and have many special properties relating to the study of the smoothness of Schubert varieties~\cite{Fan.C:B}, Kazhdan--Lusztig polynomials~\cite{Billey.S;Jones.B:A,Green.R:K}, and the decomposition of a Coxeter group into cells~\cite{Green.R;Losonczy.J:B,Shi.J:B}.

Let $w \in \FC(W)$.  Suppose that $w$ has a reduced expression beginning with $s\in S$.  Then we say that $w$ is left star reducible by $s$ with respect to $t$ to the shorter element $sw$ provided that there exists $t \in S$ such that $s$ and $t$ do not commute and $sw$ has a reduced expression beginning with $t$~\cite{Green.R:P}.  We make an analogous definition for right star reducible.  The definition of star reducible is related to C.K.~Fan's definition of cancellable in~\cite[\textsection 4]{Fan.C:A} and is also a special case of D.~Kazhdan and G.~Lusztig's notion of a star operation, which is defined for arbitrary Coxeter systems in~\cite[\textsection 10.2]{Lusztig.G:A}.

We say that $W$ is star reducible if for every fully commutative $w$, there exists a sequence $u = w_0, w_1, \ldots, w_k = w$ such that each $w_{i+1}$ is left or right star reducible to $w_i$ and $u$ is equal to a product of commuting generators.  It turns out that a Coxeter group of type $\C$ is star reducible if and only if there is an even number of generators.  However, Coxeter groups of type $B$ are star reducible regardless of the parity of the generating set~\cite[Theorem 6.3]{Green.R:P}.  In a star reducible Coxeter group, products of commuting generators form the set of fully commutative elements that are irreducible under star reductions.

In this paper, we weaken the notion of star reducible and define the non-cancellable elements, which include products of commuting generators.  Let $w \in \FC(W)$ and suppose that $w$ has a reduced expression beginning with $s\in S$.  We say that $w$ is left weak star reducible if (i) $w$ is star reducible by $s$ with respect to $t$,  and (ii) $tw$ is no longer fully commutative.  We make an analogous definition for right weak star reducible and define an element $w \in \FC(W)$ to be non-cancellable (or weak star irreducible) if it is neither left or right weak star reducible.

The non-cancellable elements of a Coxeter group $W$ are intimately related to the two-sided cells of the generalized Temperley--Lieb algebra (in the sense of Graham ~\cite{Graham.J:A}) associated to $W$.  The connection between the non-cancellable elements and the two-sided cells has been examined for types $E$ and $\widetilde{A}$ in~\cite{Fan.C:A} and~\cite{Fan.C;Green.R:A}, respectively.  This idea is also briefly touched upon for types $B$, $F$, and $H$ in~\cite{Fan.C:A}.  Due to length considerations, we will not elaborate on the connection between the non-cancellable elements and the two-sided cells.

Our motivation for studying the non-cancellable elements stems from the fact that computation involving the monomial basis elements of the generalized Tem\-per\-ley--Lieb algebra of $W$ that are indexed by non-cancellable elements is ``well-behaved.''  In fact, our classification of the non-cancellable elements in a Coxeter group of type $\C$ (Theorem~\ref{thm:affineCwsrm}) is a key component in the proof that establishes the faithfulness of a diagrammatic representation of the generalized Temperley--Lieb algebra of type $\C$, which is the focus of subsequent papers by the author.

In Section~\ref{sec:preliminaries} of this paper, we establish our notation and introduce all of the necessary terminology.  In Section~\ref{sec:type_I_and_II}, we explore some of the combinatorics of Coxeter groups of types $B$ and $\C$ and introduce the type I and type II elements, which play a central role in this paper. 

Section~\ref{sec:Bwsrm} is concerned with classifying the non-cancellable elements in a Coxeter group of type $B$ (Theorem~\ref{thm:Bwsrm}), which verifies Fan's unproved claim in~\cite[\textsection 7.1]{Fan.C:A} about the set of fully commutative elements in a Coxeter group of type $B$ having no generator appearing in the left or right descent set that can be left or right cancelled, respectively.

Using the classification of the type $B$ non-cancellable elements, we prove the main result of this paper (Theorem~\ref{thm:affineCwsrm}), which classifies the infinitely many non-cancellable elements in a Coxeter group of type $\C$.  The proof of Theorem~\ref{thm:affineCwsrm}, as well as the preparatory lemmas, rely heavily on the notation of heaps that we develop in Section~\ref{subsec:heaps}.

Lastly, in Section~\ref{sec:closing}, we expand on our discussion of our motivation for classifying the non-cancellable elements in a Coxeter group of type $\C$ and briefly discuss future research.

This paper is an adaptation of chapters 1--5 of the author's 2008 PhD thesis, titled \textit{A diagrammatic representation of an affine $C$ Temperley--Lieb algebra}~\cite{Ernst.D:A}, which was directed by Richard M. Green at the University of Colorado at Boulder.  However, some of the results presented here, especially those in Section~\ref{sec:Bwsrm}, have new and streamlined arguments.

\end{section}

%%%%%%%%%%%%%%%%%%%%%

\begin{section}{Preliminaries}\label{sec:preliminaries}

%%%%%%%%%%%%%%%%%%%%%

\begin{subsection}{Coxeter groups}\label{subsec:Coxeter_groups}

A \emph{Coxeter system} is pair $(W,S)$ consisting of a distinguished (finite) set $S$ of generating involutions and a group $W$, called a \emph{Coxeter group}, with presentation
\[
W = \langle S :(st)^{m(s, t)} = 1 \text{ for } m(s, t) < \infty \rangle,
\] 
where $m(s, s) = 1$ and $m(s, t) = m(t, s)$.  It turns out that the elements of $S$ are distinct as group elements, and that $m(s, t)$ is the order of $st$.  Given a Coxeter system $(W,S)$, the associated \emph{Coxeter graph} is the graph $\Gamma$ with vertex set $S$ and edges $\{s,t\}$ labeled with $m(s,t)$ for all $m(s,t)\geq 3$.  If $m(s,t)=3$, it is customary to leave the corresponding edge unlabeled.  Given a Coxeter graph $\Gamma$, we can uniquely reconstruct the corresponding Coxeter system $(W,S)$.  In this case, we say that the corresponding Coxeter system is of type $\Gamma$, and denote the Coxeter group and distinguished generating set by $W(\Gamma)$ and $S(\Gamma)$, respectively.

Given a Coxeter system $(W,S)$, an \emph{expression} is any product of generators from $S$.  The \emph{length} $l(w)$ of an element $w \in W$ is the minimum number of generators appearing in any expression for the element $w$.  Such a minimum length expression is called a \emph{reduced expression}.  (Any two reduced expressions for $w \in W$ have the same length.)  A product $w_{1}w_{2}\cdots w_{r}$ with $w_{i} \in W$ is called \emph{reduced} if $l(w_{1}w_{2}\cdots w_{r})=\sum l(w_{i})$.  Each element $w \in W$ can have several different reduced expressions that represent it.  Given $w \in W$, if we wish to emphasize a fixed, possibly reduced, expression for $w$, we represent it in \textsf{sans serif} font, say $\w=s_{x_1}s_{x_2}\cdots s_{x_m}$, where each $s_{x_i} \in S$.

Matsumoto's Theorem~\cite[Theorem 1.2.2]{Geck.M;Pfeiffer.G:A} says that if $w \in W$, then every reduced expression for $w$ can be obtained from any other by applying a sequence of \emph{braid moves} of the form 
\[
{\underbrace{sts \cdots }_{m(s,t)} } \mapsto {\underbrace{tst \cdots}_{m(s,t)}}
\]
where $s,t \in S$, and each factor in the move has $m(s,t)$ letters.  The \emph{support} of an element $w \in W$, denoted $\supp(w)$, is the set of all generators appearing in any reduced expression for $w$, which is well-defined by Matsumoto's Theorem.  If $\supp(w)=S$, then we say that $w$ has \emph{full support}.

Given a reduced expression $\w$ for $w \in W$, we define a \emph{subexpression} of $\w$ to be any expression obtained by deleting some subsequence of generators in the expression for $\w$.  We will refer to a consecutive subexpression of $\w$ as a \emph{subword}.

Let $w \in W$.  We write
\[
\L(w)=\{s \in S: l(sw) < l(w)\}
\]
and
\[
\R(w)=\{s \in S: l(ws) < l(w)\}.
\]
The set $\L(w)$ (respectively, $\R(w)$) is called the \emph{left} (respectively, \emph{right}) \emph{descent set} of $w$.  It turns out that $s \in \L(w)$ (respectively, $\R(w)$) if and only if $w$ has a reduced expression beginning (respectively, ending) with $s$.

The main focus of this paper will be the Coxeter systems of types $B_n$ and $\C_n$, which are defined by the following Coxeter graphs, where $n\geq 2$.

\begin{itemize}

\item[$B_n$] \begin{tabular}{l}\includegraphics[scale=.85]{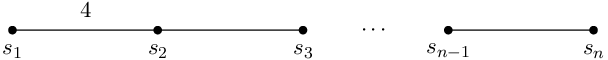}\end{tabular}

\medskip

\item[$\C_n$] \begin{tabular}{l}\includegraphics[scale=.85]{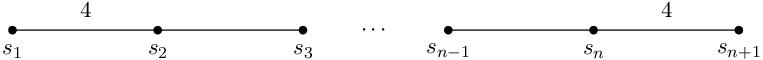}\end{tabular}

\end{itemize}

We can obtain $W(B_{n})$ from $W(\C_{n})$ by removing the generator $s_{n+1}$ and the corresponding relations~\cite[Chapter 5]{Humphreys.J:A}.  We also obtain a Coxeter group of type $B$ if we remove the generator $s_{1}$ and the corresponding relations.  To distinguish these two cases, we let $W(B_{n})$ denote the subgroup of $W(\C_{n})$ generated by $\{s_{1}, s_{2}, \dots, s_{n}\}$ and we let $W(B'_{n})$ denote the subgroup of $W(\C_{n})$ generated by $\{s_{2}, s_{3}, \dots, s_{n+1}\}$.  It is well-known that $W(\C_{n})$ is an infinite Coxeter group while $W(B_{n})$ and $W(B'_{n})$ are both finite~\cite[Chapters 2 and 6]{Humphreys.J:A}.

\end{subsection}

%%%%%%%%%%%%%%%%%%%%%

\begin{subsection}{Fully commutative elements}\label{subsec:FC}

Let $(W,S)$ be a Coxeter system of type $\Gamma$ and let $w \in W$. Following Stembridge~\cite{Stembridge.J:B}, we define a relation $\sim$ on the set of reduced expressions for $w$.  Let $\w$ and $\w'$ be two reduced expressions for $w$.  We define $\w \sim \w'$ if we can obtain $\w'$ from $\w$ by applying a single commutation move of the form $st \mapsto ts$, where $m(s,t)=2$.  Now, define the equivalence relation $\approx$ by taking the reflexive transitive closure of $\sim$.  Each equivalence class under $\approx$ is called a \emph{commutation class}. If $w$ has a single commutation class, then we say that $w$ is \emph{fully commutative}.  By Matsumoto's Theorem, an element $w$ is fully commutative if and only if no reduced expression for $w$ contains a subword of the form $sts \cdots$ of length $m(s,t) \geq 3$.  The set of all fully commutative elements of $W$ is denoted by $\FC(W)$ or $\FC(\Gamma)$.

\begin{rem}\label{rem:illegal_convex_chains}
The elements of $\FC(\C_{n})$ are precisely those whose reduced expressions avoid consecutive subwords of the following types:
\begin{enumerate}
\item $s_{i}s_{j}s_{i}$ for $|i-j|=1$ and $1< i,j < n+1$;
\item $s_{i}s_{j}s_{i}s_{j}$ for $\{i,j\}=\{1,2\}$ or $\{n,n+1\}$.
\end{enumerate}
The fully commutative elements of $W(B_{n})$ and $W(B'_{n})$ avoid the respective subwords above.
\end{rem}

In~\cite{Stembridge.J:B}, Stembridge classified the Coxeter groups that contain a finite number of fully commutative elements.  According to~\cite[Theorem 5.1]{Stembridge.J:B}, $W(\C_{n})$ contains an infinite number of fully commutative elements, while $W(B_{n})$ (and hence $W(B'_n)$) contains finitely many.  There are examples of infinite Coxeter groups that contain a finite number of fully commutative elements.  For example, Coxeter groups of type $E_n$ for $n\geq 9$ are infinite, but contain only finitely many fully commutative elements~\cite[Theorem 5.1]{Stembridge.J:B}.

\end{subsection}

%%%%%%%%%%%%%%%%%%%%%

\begin{subsection}{Non-cancellable elements}\label{subsec:non-cancellable}

The notion of a star operation was originally defined by Kazhdan and Lusztig in~\cite[\textsection 4.1]{Kazhdan.D;Lusztig.G:A} for simply laced Coxeter systems (i.e., $m(s,t)\leq3$ for all $s,t \in S$) and was later generalized to arbitrary Coxeter systems in~\cite[\textsection 10.2]{Lusztig.G:A}.  If $I=\{s,t\}$ is a pair of noncommuting generators for $W$, then $I$ induces four partially defined maps from $W$ to itself, known as star operations. A star operation, when it is defined, respects the partition $W = \FC(W)\ \dot{\cup}\  (W \setminus \FC(W) )$ of the Coxeter group, and increases or decreases the length of the element to which it is applied by 1.  For our purposes, it is enough to define star operations that decrease length by 1, and so we will not develop the full generality.

Suppose that $(W,S)$ is an arbitrary Coxeter system of type $\Gamma$.  Let $w \in W$ and suppose that $s \in \L(w)$.  We define $w$ to be \emph{left star reducible by $s$ with respect to $t$} to $sw$ if there exists $t \in \L(sw)$ with $m(s,t) \geq 3$.  We analogously define \emph{right star reducible by $s$ with respect to $t$}.  Observe that if $m(s,t)\geq 3$, then $w$ is left (respectively, right) star reducible by $s$ with respect to $t$ if and only if $w=stv$ (respectively, $w=vts$), where the product is reduced.  We say that $w$ is \emph{star reducible} if it is either left or right star reducible by some $s\in S$.

We now introduce the concept of weak star reducible, which is related to Fan's notion of cancellable in~\cite{Fan.C:A}.   If $w \in \FC(W)$, then $w$ is \emph{left weak star reducible by $s$ with respect to $t$} to $sw$ if (i) $w$ is left star reducible by $s$ with respect to $t$, and (ii) $tw \notin \FC(W)$.  Observe that (i) implies that $m(s,t) \geq 3$ and that $s \in \L(w)$.  Furthermore, (ii) implies that $l(tw)>l(w)$.  Also, note that we are restricting our definition of weak star reducible to the set of fully commutative elements.  We analogously define \emph{right weak star reducible by $s$ with respect to $t$}.  If $w$ is either left or right weak star reducible by some $s\in S$, we say that $w$ is \emph{weak star reducible}.  Otherwise, we say that $w\in \FC(W)$ is \emph{non-cancellable} or \emph{weak star irreducible} (or simply \emph{irreducible}).

\begin{ex}
Consider $w,w' \in \FC(\C_{n})$ having reduced expressions $\w=s_{1}s_{2}s_{1}$ and $\w'=s_{1}s_{2}$, respectively.  We see that $w$ is left (respectively, right) weak star reducible by $s_{1}$ with respect to $s_{2}$ to $s_2 s_1$ (respectively, $s_1 s_2$), and so $w$ is not non-cancellable.  However, $w'$ is non-cancellable.
\end{ex}    

\begin{rem}
We make a few observations regarding weak star operations.
\begin{enumerate}

\item If $w\in \FC(W)$ and $s \in \L(w)$ (respectively, $\R(w)$), it is clear that $sw$ (respectively, $ws$) is still fully commutative.  This implies that if $w \in \FC(W)$ is left or right weak star reducible to $u$, then $u$ is also fully commutative.  

\item It follows immediately from the definition that if $w$ is weak star reducible to $u$, then $w$ is also star reducible to $u$.  However, there are examples of fully commutative elements that are star reducible, but not weak star reducible.  For example, consider $w=s_{1}s_{2} \in \FC(B_{2})$.  We see that $w$ is star reducible, but not weak star reducible since $tw$ and $wt$ are still fully commutative for any $t \in S$.  However, observe that in simply laced Coxeter systems (i.e., $m(s,t)\leq 3$ for all $s,t\in S$), star reducible and weak star reducible are equivalent.

\item If $w \in \FC(\C_{n})$, then $w$ is left weak star reducible by $s$ with respect to $t$ if and only if $w=stv$ (reduced) when $m(s,t)=3$, or $w=stsv$ (reduced) when $m(s,t)=4$.   Again, observe that the characterization above applies to $\FC(B_{n})$ and $\FC(B'_{n})$. 

\end{enumerate}

\end{rem}

\end{subsection}

%%%%%%%%%%%%%%%%%%%%%

\begin{subsection}{Heaps}\label{subsec:heaps}

Every reduced expression can be associated with a partially ordered set called a heap that will allow us to visualize a reduced expression  while preserving the essential information about the relations among the generators.  The theory of heaps was introduced in~\cite{Viennot.G:A} by Viennot and visually captures the combinatorial structure of the Cartier--Foata monoid of~\cite{Cartier.P;Foata.D:A}.  In~\cite{Stembridge.J:B} and~\cite{Stembridge.J:A}, Stembridge studied heaps in the context of fully commutative elements, which is our motivation here.  

Although heaps will be useful for visualizing the arguments throughout the remainder of this paper, we will not exploit their full utility until Section~\ref{sec:Cwsrm}, where we classify the non-cancellable elements of type $\C$.  In this section, we mimic the development found in~\cite{Billey.S;Jones.B:A},~\cite{Billey.S;Warrington.G:A}, and~\cite{Stembridge.J:B}.  

Let $(W,S)$ be a Coxeter system.  Suppose $\w = s_{x_1} \cdots s_{x_r}$ is a fixed reduced expression for $w \in W$.  As in~\cite{Stembridge.J:B}, we define a partial ordering on the indices $\{1, \dots, r\}$ by the transitive closure of the relation $\lessdot$ defined via $j \lessdot i$ if $i < j$ and $s_{x_i}$ and $s_{x_j}$ do not commute.  In particular, $j \lessdot i$ if $i < j$ and $s_{x_i} = s_{x_j}$ (since we took the transitive closure).  This partial order is referred to as the \emph{heap} of $\w$, where $i$ is labeled by $s_{x_i}$.  It follows from~\cite[Proposition 2.2]{Stembridge.J:B} that heaps are well-defined up to commutativity class.  That is, if $\w$ and $\w'$ are two reduced expressions for $w \in W$ that are in the same commutativity class, then the labeled heaps of $\w$ and $\w'$ are equal.  In particular, if $w$ is fully commutative, then it has a single commutativity class, and so there is a unique heap associated to $w$.

\begin{ex}\label{ex:first.heap}
Let $\w = s_3 s_2 s_1 s_2 s_5s_{4}s_{6}s_{5}$ be a reduced expression for $w \in \FC(\C_{5})$.  We see that $\w$ is indexed by $\{1, 2, 3, 4, 5, 6, 7, 8\}$.  As an example, $3 \lessdot 2$ since $2 < 3$ and the second and third generators do not commute.  The labeled Hasse diagram for the unique heap poset of $w$ is shown below.
\[
\includegraphics[scale=.85]{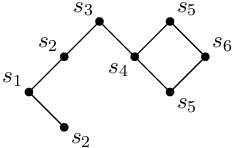}
\]
\end{ex}

Let $\w$ be a fixed reduced expression for $w \in W(\C_{n})$.  As in~\cite{Billey.S;Jones.B:A} and~\cite{Billey.S;Warrington.G:A}, we will represent a heap for $\w$ as a set of lattice points embedded in $\{1,2,\ldots,n+1\} \times \mathbb{N}$.  To do so, we assign coordinates (not unique) $(x,y) \in \{1,2,\ldots, n+1\} \times \mathbb{N}$ to each entry of the labeled Hasse diagram for the heap of $\w$ in such a way that:
\begin{enumerate}
\item An entry with coordinates $(x,y)$ is labeled $s_i$ in the heap if and only if $x = i$; 

\item An entry with coordinates $(x,y)$ is greater than an entry with coordinates $(x',y')$ in the heap if and only if $y > y'$.
\end{enumerate}

Recall that a finite poset is determined by its covering relations.  In the case of $\C_{n}$ (and any straight line Coxeter graph), it follows from the definition that $(x,y)$ covers $(x',y')$ in the heap if and only if $x = x' \pm 1$, $y > y'$, and there are no entries $(x'', y'')$ such that $x'' \in \{x, x'\}$ and $y'< y'' < y$.  This implies that we can completely reconstruct the edges of the Hasse diagram and the corresponding heap poset from a lattice point representation. The lattice point representation of a heap allows us to visualize potentially cumbersome arguments.  Note that our heaps are upside-down versions of the heaps that appear in in~\cite{Billey.S;Jones.B:A} and~\cite{Billey.S;Warrington.G:A} and several other papers.  That is, in this paper entries on top of a heap correspond to generators occurring to the left, as opposed to the right, in the corresponding reduced expression.  However, our convention aligns more naturally with the typical conventions of diagram algebras that are motivating the results of this paper.
%There's a typo above: "in in"

Let $\w$ be a reduced expression for $w \in W(\C_{n})$.  We let $H(\w)$ denote a lattice representation of the heap poset in $\{1,2,\ldots,n+1\} \times \N$ described in the preceding paragraph.  If $w$ is fully commutative, then the choice of reduced expression for $w$ is irrelevant, in which case, we will often write $H(w)$ (note the absence of \textsf{sans serif} font) and we will refer to $H(w)$ as the heap of $w$.

Given a heap, there are many possible coordinate assignments, yet the $x$-co\-or\-di\-nates for each entry will be fixed for all of them.  In particular, two entries labeled by the same generator may only differ by the amount of vertical space between them while maintaining their relative vertical position to adjacent entries in the heap.

Let $\w=s_{x_1}\cdots s_{x_r}$ be a reduced expression for $w \in \FC(\C_{n})$.  If $s_{x_i}$ and $s_{x_j}$ are adjacent generators in the Coxeter graph with $i<j$, then we must place the point labeled by $s_{x_i}$ at a level that is \textit{above} the level of the point labeled by $s_{x_j}$.  Because generators that are not adjacent in the Coxeter graph do commute, points whose $x$-coordinates differ by more than one can slide past each other or land at the same level.  To emphasize the covering relations of the lattice representation we will enclose each entry of the heap in a rectangle in such a way that if one entry covers another, the rectangles overlap halfway.

\begin{ex}\label{second.heap.ex}
Let $w$ be as in Example~\ref{ex:first.heap}.  Then one possible representation for $H(w)$ is as follows.
\[
\includegraphics[scale=.85]{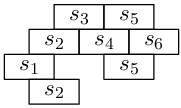}
\]
\end{ex}

When $w$ is fully commutative, we wish to make a canonical choice for the representation $H(w)$ by assembling the entries in a particular way.  To do this, we give all entries corresponding to elements in $\L(w)$ the same vertical position and all other entries in the heap should have vertical position as high as possible.  Note that our canonical representation of heaps of fully commutative elements corresponds precisely to the unique heap factorization of~\cite[Lemma 2.9]{Viennot.G:A} and to the Cartier--Foata normal form for monomials~\cite{Cartier.P;Foata.D:A,Green.R:P}.  In Example~\ref{second.heap.ex}, the representation of $H(w)$ that we provided is the canonical representation.  When illustrating heaps, we will adhere to this canonical choice, and when we consider the heaps of arbitrary reduced expressions, we will only allude to the relative vertical positions of the entries, and never their absolute coordinates.  

Given a canonical representation of a heap, it makes sense to refer to the $k$th row of the heap, and we will do this when no confusion will arise.  Note that for fully commutative elements, the first row of the heap corresponds to the left descent set.  If $w \in \FC(\C_{n})$, let $\r_{k}$ denote the $k$th row of the canonical representation for $H(w)$.  We will write $s_{i} \in \r_{k}$ to mean that there is an entry occurring in the $k$th row labeled by $s_{i}$.  If $\r_{k}$ consists entirely of entries labeled by $s_{x_{1}}, s_{x_{2}}, \dots, s_{x_{m}}$, then we will write $\r_{k}=s_{x_{1}} s_{x_{2}} \cdots s_{x_{m}}$.  

Let $w \in \FC(\C_n)$ have reduced expression $\w=s_{x_1}\cdots s_{x_r}$ and suppose $s_{x_i}$ and $s_{x_j}$ equal the same generator $s_k$, so that the corresponding entries have $x$-coordinate $k$ in $H(w)$.  We say that $s_{x_i}$ and $s_{x_j}$ are \emph{consecutive} if there is no other occurrence of $s_{k}$ occurring between them in $\w$.  In this case, $s_{x_i}$ and $s_{x_j}$ are consecutive in $H(w)$, as well.  

Let $\w=s_{x_{1}} \cdots s_{x_{r}}$ be a reduced expression for $w \in W(\C_{n})$.  We define a heap $H'$ to be a \emph{subheap} of $H(\w)$ if $H'=H(\w')$, where $\w'=s_{y_1}s_{y_2} \cdots s_{y_k}$ is a subexpression of $\w$.  We emphasize that the subexpression need not be a subword (i.e., a consecutive subexpression).  

Recall that a subposet $Q$ of $P$ is called convex if $y \in Q$ whenever $x < y < z$ in $P$ and $x, z \in Q$.  We will refer to a subheap as a \emph{convex subheap} if the underlying subposet is convex.  

\begin{ex}
As an example, let $\w= s_3 s_2 s_1 s_2 s_5s_{4}s_{6}s_{5}$ as in Example~\ref{ex:first.heap}.  Now, let $\w'=s_{5}s_{4}s_{5}$ be the subexpression of $\w$ that results from deleting all but fifth, sixth, and last generators of $\w$.  Then $H(\w')$ equals 
\[
\includegraphics[scale=.85]{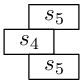}
\]
and is a subheap of $H(\w)$, but is not convex since there is an entry in $H(\w)$ labeled by $s_{6}$ occurring between the two consecutive occurrences of $s_{5}$ that does not occur in $H(\w')$.  However, if we do include the entry labeled by $s_{6}$, then
\[
\includegraphics[scale=.85]{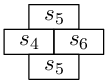}
\]
is a convex subheap of $H(\w)$.
\end{ex}

From this point on, if there can be no confusion, we will not specify the exact subexpression that a subheap arises from.

The following fact is implicit in the literature (in particular, see the proof of~\cite[Proposition 3.3]{Stembridge.J:B}) and follows easily from the definitions.

\begin{prop}
Let $w \in \FC(W)$.  Then $H'$ is a convex subheap of $H(w)$ if and only if $H'$ is the heap for some subword of some reduced expression for $w$.   \qed
\end{prop}

It will be extremely useful for us to be able to recognize when a heap corresponds to a fully commutative element in $W(\C_{n})$.  The following lemma follows immediately from Remark~\ref{rem:illegal_convex_chains} and is also a special case of~\cite[Proposition 3.3]{Stembridge.J:B}.

\begin{lem}\label{lem:impermissible.heap.configs}
Let $w \in \FC(\C_{n})$.  Then $H(w)$ cannot contain any of the following convex subheaps:
\begin{center}
\begin{tabular}[c]{c}
\includegraphics[scale=.85]{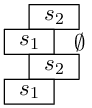}
\end{tabular},
\begin{tabular}[c]{c}
\includegraphics[scale=.85]{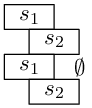}
\end{tabular},
\begin{tabular}[c]{c}
\includegraphics[scale=.85]{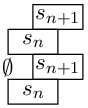}
\end{tabular},
\begin{tabular}[c]{c}
\includegraphics[scale=.85]{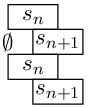}
\end{tabular},
\begin{tabular}[c]{c}
\includegraphics[scale=.85]{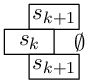}
\end{tabular},
\begin{tabular}[c]{c}
\includegraphics[scale=.85]{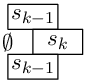}
\end{tabular},
\end{center}
where $1<k<n+1$ and we use $\emptyset$ to emphasize that no element of the heap occupies the corresponding position.  \qed
\end{lem}

We conclude this section with an observation regarding heaps and weak star reductions.  Let $\w=s_{x_1}\cdots s_{x_r}$ be a reduced expression for $w \in \FC(\C_{n})$.  Then $w$ is left weak star reducible by $s$ with respect to $t$ if and only if
\begin{enumerate}
\item there is an entry in $H(\w)$ labeled by $s$ that is not covered by any other entry; and

\item the heap $H(t\w)$ contains one of the convex subheaps of Lemma~\ref{lem:impermissible.heap.configs}.
\end{enumerate}
Of course, we have an analogous statement for right weak star reducible.

\end{subsection}

\end{section}

%%%%%%%%%%%%%%%%%%%%%

\begin{section}{The type {I} and type II elements of a Coxeter group of type $\C$}\label{sec:type_I_and_II}

In this section, we explore some of the combinatorics of Coxeter groups of types $B$ and $\C$.  Our immediate goal is to define two classes of fully commutative elements of $W(\C_{n})$ that play a central role in the remainder of this paper.  Most of these elements will turn out to be on our list of non-cancellable elements appearing in Section~\ref{sec:Cwsrm}.

%%%%%%%%%%%%%%%%%%%%%

\begin{subsection}{The type I elements}

Let $w \in \FC(\C_{n})$.  We define $n(w)$ to be the maximum integer $k$ such that $w$ has a reduced expression of the form $w = u x v$ (reduced), where $u, x, v \in \FC(\C_{n})$, $l(x)=k$, and $x$ is a product of commuting generators.  Note that $n(w)$ may be greater than the size of any row in the canonical representation of $H(w)$.  Also, it is known that $n(w)$ is equal to the size of a maximal antichain in the heap poset for $w$~\cite[Lemma 2.9]{Shi.J:C}.

\begin{defn}\label{def:zigzags}
Define the following elements of $W(\C_{n})$.
\begin{enumerate}
\item If $i<j$, let
\[
\z_{i,j}=s_{i}s_{i+1}\cdots s_{j-1}s_{j}
\]
and
\[
\z_{j,i}=s_{j}s_{j-1}\cdots s_{i-1}s_{i}.
\]
We also let $\z_{i,i}=s_{i}$.

\item If $1< i \leq n+1$ and $1 < j \leq n+1$, let
\[
\z^{L,2k}_{i,j}=\z_{i,2}(\z_{1,n}\z_{n+1,2})^{k-1}\z_{1,n}\z_{n+1,j}.
\]
\item If $1< i \leq n+1$ and $1 \leq j < n+1$, let
\[
\z^{L,2k+1}_{i,j}=\z_{i,2}(\z_{1,n}\z_{n+1,2})^{k}\z_{1,j}.
\]

\item If $1\leq i < n+1$ and $1 \leq j <  n+1$, let
\[
\z^{R,2k}_{i,j}=\z_{i,n}(\z_{n+1,2}\z_{1,n})^{k-1}\z_{n+1,2}\z_{1,j}.
\]
	
\item If $1\leq i < n+1$ and $1 < j \leq  n+1$, let 
\[
\z^{R,2k+1}_{i,j}=\z_{i,n}(\z_{n+1,2}\z_{1,n})^{k}\z_{n+1,j}.
\]

\end{enumerate}
If $w \in W(\C_n)$ is equal to one of the elements in (1)--(5), then we say that $w$ is of \emph{type I}.
\end{defn}

The notation for the type I elements looks more cumbersome than the underlying concept.  The notation is motivated by the zigzagging shape of the corresponding heaps.  The index $i$ tells us where to start and the index $j$ tells us where to stop.  The L (respectively, R) tells us to start zigzagging to the left (respectively, right).  Also, $2k+1$ (respectively, $2k$) indicates the number of times we should encounter an end generator (i.e., $s_{1}$ or $s_{n+1}$) after the first occurrence of $s_{i}$ as we zigzag through the generators.  If $s_{i}$ is an end generator, it is not included in this count.  However, if $s_{j}$ is an end generator, it is included.  

\begin{ex}
If $1<i,j\leq n+1$, then
\[
H\(\z^{L,2k}_{i,j}\)=\begin{tabular}[c]{c}\includegraphics[scale=.85]{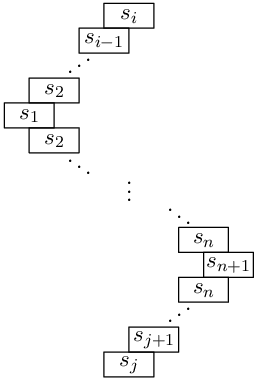}\end{tabular},
\]
where we encounter an entry labeled by either $s_{1}$ or $s_{n+1}$ a combined total of $2k$ times if $i\neq n+1$ and $2k+1$ times if $i=n+1$.  
\end{ex}

Every type I element is rigid, in the sense that each has a unique reduced expression.  This implies that every type I element is fully commutative (there are no relations of any kind to apply).  Furthermore, it is clear from looking at the heaps for the type I elements that if $w$ is of type I, then $n(w)=1$.  Conversely, it follows by induction on $l(w)$ that if $n(w)=1$ for some $w \in \FC(\C_{n})$, then $w$ must be of type I.  Lastly, note that there are an infinite number of type I elements since there is no limit to the zigzagging that their corresponding heaps can do.

The discussion in the previous paragraph verifies the following proposition.

\begin{prop}\label{prop:zigzags}
If $w \in W(\C_{n})$ is of type I, then $w$ is fully commutative with $n(w)=1$.  Conversely, if $n(w)=1$, then $w$ is one of the elements on the list in Definition~\ref{def:zigzags}.   \qed
\end{prop}

\end{subsection}

%%%%%%%%%%%%%%%%%%%%%

\begin{subsection}{The type II elements}

It will be helpful for us to define $\lambda=\lceil \frac{n-1}{2}\rceil$.  Then regardless of whether $n$ is odd or even, $2\lambda$ will always be the largest even number in $\{1,2,\dots,n+1\}$.  Similarly, $2\lambda+1$ will always be the largest odd number in $\{1,2,\dots, n+1\}$.

\begin{defn}
Define $\O=\{1,3, \dots, 2\lambda-1, 2\lambda+1\}$ and $\E=\{2, 4, \dots, 2\lambda-2, 2\lambda\}$.  (Note that $\O$ (respectively, $\E$) consists of all of the odd (respectively, even) indices amongst $\{1, 2, \dots,n+1\}$.)  Let $i$ and $j$ be of the same parity with $i<j$.  We define 
\[
\x_{i,j}=s_{i}s_{i+2}\cdots s_{j-2}s_{j}.
\]
Also, define
\[
\x_{\O}=\x_{1,2\lambda+1}=s_{1}s_{3}\cdots s_{2\lambda-1}s_{2\lambda+1},
\]
and
\[
\x_{\E}=\x_{2,2\lambda}=s_{2}s_{4}\cdots s_{2\lambda-2}s_{2\lambda}.
\]
If $w \in W(\C_n)$ is equal to a finite alternating product of $\x_{\O}$ and $\x_{\E}$, then we say that $w$ is of \emph{type II}.  (It is important to point out that the corresponding expressions are indeed reduced.)
\end{defn}

\begin{ex}
If $n$ is even, then
\[
H\(\x_{\E}(\x_{\O}\x_{\E})^{k}\)=\begin{tabular}[c]{c}
\includegraphics[scale=.85]{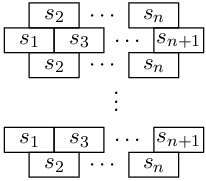}\end{tabular},
\]
where the canonical representation has $2k+1$ rows.
\end{ex}

The next proposition follows immediately since the heaps of the type II elements avoid the impermissible configurations of Lemma~\ref{lem:impermissible.heap.configs}.

\begin{prop}
Let $w\in W(\C_{n})$ be of type II. Then $w$ is fully commutative. Moreover, if $w$ is not equal to $\x_{\E}$ when $n$ is even, then $n(w)=\lambda$.\footnote{The published version of this paper appearing in \emph{Int. Electron. J. Algebra} 8, 2010 did not exclude the case when $w$ is equal to $\x_{\E}$ for even $n$. The same error is contained in~\cite{Ernst.D:A}. Thankfully, this error has no impact on the remaining results in this paper.} \qed
%If $w \in W(\C_{n})$ is of type II, then $w$ is fully commutative with $n(w)=\lambda$.   
\end{prop}

It is quickly seen by inspecting the heaps for the type II elements that if $w$ is of type II, then $w$ is non-cancellable.  Note that if $w \in \FC(\C_{n})$, then $\lambda$ is the maximum value that $n(w)$ can take.  Furthermore, there are infinitely many type II elements.  However, not every fully commutative element with $n$-value $\lambda$ is of type II.

Note that if $n$ is even, then every $\(\x_{\O}\x_{\E}\)^{k}\x_{\O}$ with $k>0$ is not star reducible.  This fact is implicit in~\cite{Green.R:P} and is easily verified.  It follows from our classification of the type $\C$ non-cancellable elements (see Theorem~\ref{thm:affineCwsrm}) that these elements are the only non-star reducible elements in $\FC(\C_{n})$ (with $n$ even) other than products of commuting generators; all other non-cancellable elements are star reducible.

\end{subsection}

\end{section}

%%%%%%%%%%%%%%%%%%%%%

\begin{section}{The type $B$ non-cancellable elements}\label{sec:Bwsrm}

The goal of this section is to classify the non-cancellable elements of $W(B_{n})$.  To accomplish this task, we shall make use of a normal form for reduced expressions in a Coxeter system of type $B$.

%%%%%%%%%%%%%%%%%%%%%

\begin{subsection}{Preparatory lemmas}

Mimicking~\cite[\textsection 2.1]{Green.R;Losonczy.J:E}, define 
\[
W^{(k)}=\{w\in W(B_{k}): 1\leq i<k \implies \ell(s_{i}w)>\ell(w)\}.
\]
Then $W^{(k)}$ is a set of minimum length right coset representatives for the subgroup $W(B_{k-1})$ of $W(B_{k})$, and $\ell(uv)=\ell(u)+\ell(v)$ for all $u\in W(B_{k-1})$ and $v\in W(B_{k})$ (see~\cite[\textsection 5.12]{Humphreys.J:A}).  It is an easy exercise to show that the elements of $W^{(k)}$ are given by
\[
\{e, s_{k},\z_{k,k-1},\z_{k,k-2},\ldots,\z_{k,1},\z_{k,2}^{L,1},\z_{k,3}^{L,1},\ldots,\z_{k,k}^{L,1}\}.
\]
(One way this can be established is by working with the signed permutation representation of $W(B_{k})$.  Also, see~\cite[\textsection 2.1]{Green.R;Losonczy.J:E}.)

\begin{lem}\label{lem:normal form}
Let $w \in W(B_{n})$.  Then $w$ has a unique reduced decomposition $w=w_{1}w_{2}\cdots w_{n}$, where each $w_{k}\in W^{(k)}$.
\end{lem}

\begin{proof}
See proof of Lemma 2.1.1 in~\cite{Green.R;Losonczy.J:E}.
\end{proof}

We will refer to the unique reduced decomposition of Lemma~\ref{lem:normal form} as the \emph{normal form factorization} for $w$.

The next two lemmas play a crucial role in the proof of Theorem~\ref{thm:Bwsrm}.

\begin{lem}\label{lem:above_diagonal}
Let $w\in \FC(B_{n})$ have normal form factorization $w=w_{1}w_{2}\cdots w_{n}$.  If there exists $k$ such that $w_{k}=\z_{k,1}$, then for each $1\leq i<k$, $w_{i}$ equals the identity $e$ or equals $\z_{i,1}$.
\end{lem}

\begin{proof}
For sake of a contradiction, assume otherwise.  Choose the largest $i$ such that $1\leq i<k$ and $w_{i}$ is not equal to either the identity or $\z_{i,1}$.  First, observe that we must have $i>1$ since $W^{(1)}=\{e,s_1\}$.   By how we chose $i$, there must exist $m$ with $i<m\leq k$ such that $w_m=\z_{m,1}$.   Choose the smallest such $m$, so that $w_{l}=e$ for all $i<l<m$.  Then the only possibilities are that $w_{i}=\z_{i,j}$ with $1<j<i$, or $w_{i}=\z_{i,j}^{L,1}$ with $1<j\leq i$.  Furthermore, $s_{j+1}s_j$ is a subword of some reduced expression for $w_m$.  This implies that some reduced expression for $w$ would contain the subword $s_{j}s_{j+1}s_{j}$, where the first occurrence of $s_j$ comes from $w_i$ while $s_{j+1}s_j$ comes from $w_m$.  This violates $w$ being fully commutative.
\end{proof}

\begin{lem}\label{lem:above_triangle}
Let $w\in \FC(B_{n})$ have normal form factorization $w=w_{1}w_{2}\cdots w_{n}$ and suppose that $w$ is non-cancellable such that $w_{n-i}=\z_{n-i,n-2i}$ for all $1\leq i \leq k$ for some $k<\frac{n-1}{2}$.  Then $w_{n-(k+1)}$ is equal to $\z_{n-(k+1),n-2(k+1)}$ or $\z_{n-(k+1),1}$.
\end{lem}

\begin{proof}
Since $s_{n-k}$ does not appear in the support of $w_{1}w_{2}\cdots w_{n-(k+1)}$ and $w$ is not left weak star reducible, it must be the case that $w_{n-(k+1)}\neq e$; otherwise, $w$ is left weak star reducible by $s_{n-k}$ with respect to $s_{n-(k-1)}$.  Since $w$ is fully commutative and not right weak star reducible, it quickly follows that the only possibilities for $w_{n-(k+1)}$ in $W^{(n-(k+1))}$ are $\z_{n-(k+1),n-2(k+1)}$ or $\z_{n-(k+1),1}$.
\end{proof}

\end{subsection}

%%%%%%%%%%%%%%%%%%%%%

\begin{subsection}{Classification of the type $B$ non-cancellable elements}

The next theorem verifies Fan's unproved claim in~\cite[\textsection 7.1]{Fan.C:A} about the set of $w \in \FC(B_{n})$ having no element of $\L(w)$ or $\R(w)$ that can be left or right cancelled, respectively.  

\begin{thm}\label{thm:Bwsrm}
Let $w \in \FC(B_{n})$.  Then $w$ is non-cancellable if and only if $w$ is equal to either a product of commuting generators, $s_{1}s_{2}u$, or $s_{2}s_{1}u$, where $u$ is a product of commuting generators with $s_{1}, s_{2}, s_{3} \notin \supp(u)$.  We have an analogous statement for $\FC(B'_{n})$, where $s_{1}$ and $s_{2}$ are replaced with $s_{n+1}$ and $s_{n}$, respectively.
\end{thm}

\begin{proof}
First, observe that if $w$ is non-cancellable in $W(B_{n'})$ for $n' < n$, then $w$ is also non-cancellable in $W(B_{n})$ when considered as an element of the larger group.  Also, we see that every element on our list is, in fact, non-cancellable.  It remains to show that our list is complete. We induct on the rank $n$.  

For the base case, consider $n=2$.  An exhaustive check verifies that the only non-cancellable elements in $W(B_{2})$ are $s_{1}$, $s_{2}$, $s_{1}s_{2}$, and $s_{2}s_{1}$, which agrees with the statement of the theorem.  

For the inductive step, assume that for all $n' \leq n-1$, our list is complete.  Let $w \in \FC(B_{n})$ and assume that $w$ is non-cancellable with normal form factorization $w=w_{1}w_{2}\cdots w_{n}$.  If $s_{n}\notin \supp(w)$, then we are done by induction.  So, assume that $s_{n}\in \supp(w)$.  In this case, $s_{n}\in\supp(w_{n})$, but $s_{n}\notin \supp(w_{i})$ for all $1\leq i <n$.  Since $w$ is not right weak star redudible, there are only two possibilities: (1) $w_n=\z_{n,1}$, or (2) $w_{n}=s_{n}$.  

Case (1): Suppose that $w_{n}=\z_{n,1}$.  Then we may apply Lemma~\ref{lem:above_diagonal} and conclude that either (a) $w_{i}=e$ for all $i<n$, or (b) $w_{i}=\z_{i,1}$ for some $i<n$.  If we are in situation (a), then $w$ would be left weak star reducible by $s_{n}$ with respect to $s_{n-1}$.  Assume that (b) occurs and choose the largest such $i$, so that $w_{i}=\z_{i,1}$ while $w_{j}=e$ for $i<j<n$.  In this case, $w$ would be right weak star reducible by $s_{1}$ with respect to $s_{2}$.  Regardless, we contradict $w$ being non-cancellable.  Thus, we must be in case (2).

Case (2):  Now, assume that $w_{n}=s_{n}$, and for sake of a contradiction, assume that $w_{n-1}\neq e$.  This implies that
\[
w_{n-1}\in \{s_{n-1},\z_{n-1,n-2},\ldots,\z_{n-1,1},\z_{n-1,2}^{L,1},\ldots,\z_{n-1,n-1}^{L,1}\}.
\]
If $w_{n-1}=s_{n-1}$, then $w$ would be right weak star reducible by $s_{n}$ with respect to $s_{n-1}$.  If $w_{n-1}=\z_{n-1,k}$ for $k\in \{2,\ldots,n-3\}$, then $w$ would be right weak star reducible by $s_{k}$ with respect to $s_{k+1}$.  Similarly, if $w_{n-1}=\z_{n-1,k}^{L,1}$ with $k\in \{2,\ldots,n-2\}$, then $w$ would be right weak star reducible by $s_{k}$ with respect to $s_{k-1}$.  Also, if $w_{n-1}=\z_{n-1,n-1}^{L,1}$, then $w$ is right weak star reducible by $s_{n}$ with respect to $s_{n-1}$.  The only remaining possibilities are: (a) $w_{n-1}=\z_{n-1,n-2}$, or (b) $w_{n-1}=\z_{n-1,1}$.

(a)  Suppose that $w_{n-1}=\z_{n-1,n-2}$.  By making repeated applications of Lemma~\ref{lem:above_triangle}, we can conclude that there exists $k$ such that $w_{k}=\z_{k,1}$.   Choose the largest such $k$.  Then $s_{k+1}\in \L(w_{k+1})$.  By Lemma~\ref{lem:above_diagonal}, $w_{k-1}$ is equal to either $e$ or $\z_{k-1,1}$.  If $w_{k-1}=e$, then $w$ would be left weak star reducible by $s_{k}$ with respect to $s_{k+1}$.  Yet, if $w_{k-1}=\z_{k-1,1}$, we would have $w$ right weak star reducible by $s_{1}$ with respect to $s_{2}$.  In either case, we contradict $w$ being non-cancellable.

(b) Lastly, assume that $w_{n-1}=\z_{n-1,1}$.  In this case, we can apply Lemma~\ref{lem:above_diagonal} and conclude that $w_{n-2}$ is equal to either $e$ or $\z_{n-2,1}$.  If $w_{n-2}=e$, then $w$ would be left weak star reducible by $s_{n-1}$ with respect to $s_{n}$.  On the other hand, if $w_{n-2}=\z_{n-2,1}$, then $w$ is right weak star reducible by $s_{1}$ with respect to $s_{2}$.  Again, we contradict $w$ being non-cancellable.

Therefore, it must be the case that $w_{n-1}=e$, which implies that $s_{n-1}\notin \supp(w)$.  In this case, we can apply the induction hypothesis to $w_{1}w_{2}\cdots w_{n-1}$ and conclude that $w$ is one of the elements on our list (since $s_{n}$ would commute with all the elements in $\supp(w_{1}w_{2}\cdots w_{n-1})$).
\end{proof}

\end{subsection}

\end{section}

%%%%%%%%%%%%%%%%%%%%%

\begin{section}{The type $\C$ non-cancellable elements}\label{sec:Cwsrm}

In this section, we will classify the non-cancellable elements of $W(\C_{n})$.  

%%%%%%%%%%%%%%%%%%%%%

\begin{subsection}{Statement of theorem}

The following theorem is the main result of this paper.

\begin{thm}\label{thm:affineCwsrm}
Let $w \in \FC(\C_{n})$.  Then $w$ is non-cancellable if and only if $w$ is equal to one of the elements on the following list.
\begin{enumerate}[label=\rm{(\roman*)}]
\item $uv$, where $u$ is a type $B$ non-cancellable element and $v$ is a type $B'$ non-cancellable element with $\supp(u)\cap \supp(v)=\emptyset$;
\item $\z^{R,2k}_{1,1}$, $\z^{L,2k}_{n+1,n+1}$, $\z^{L,2k+1}_{n+1,1}$, and $\z^{R,2k+1}_{1,n+1}$; 
\item any type II element.
\end{enumerate}
\end{thm}

The elements listed in (i) include all possible products of commuting generators.  This includes $\x_{\O}$ and $\x_{\E}$, which are also included in (iii).  The elements listed in (ii) are the type I elements having left and right descent sets equal to one of the end generators.

\end{subsection}

%%%%%%%%%%%%%%%%%%%%%

\begin{subsection}{More preparatory lemmas}

The proof of Theorem~\ref{thm:affineCwsrm} requires several technical lemmas whose proofs rely heavily on the heap notation that we developed in Section~\ref{subsec:heaps}.

Before proceeding, we make a comment on notation.  When representing convex subheaps of $H(w)$ for $w\in \FC(\C_n)$, we will use the symbol $\emptyset$ to emphasize the absence of an entry in this location of $H(w)$.  It is important to note that the occurrence of the symbol $\emptyset$  implies that an entry from the canonical representation of $H(w)$ cannot be shifted vertically from above or below to occupy the location of the symbol $\emptyset$.  If we enclose a region by a dotted line and label the region with $\emptyset$, we are indicating that no entry of the heap may occupy this region.

We will make frequent use of the following lemma, which allows us to determine whether an element is of type I.

\begin{lem}\label{lem:zigzag}
Let $w \in \FC(\C_{n})$.  Suppose that $w$ has a reduced expression having one of the following fully commutative elements as a subword:  
\begin{enumerate}[label=\rm{(\roman*)}]
\item $\z^{L,2}_{2,n}=s_{2}s_{1}s_{2}s_{3} \cdots s_{n-1}s_{n}s_{n+1}s_{n}$,
\item $\z^{R,2}_{n,2}=s_{n}s_{n+1}s_{n}s_{n-1}\cdots s_{3}s_{2}s_{1}s_{2}$,
\item $\z^{R,2}_{1,1}=s_{1}s_{2}\cdots s_{n}s_{n+1}s_{n}\cdots s_{2}s_{1}$,
\item $\z^{L,2}_{n+1,n+1}=s_{n+1}s_{n}\cdots s_{2}s_{1}s_{2}\cdots s_{n}s_{n+1}$ .
\end{enumerate}
Then $w$ is of type I.
\end{lem}

\begin{proof}
One quickly sees that if $w$ has any of the above reduced expressions as a subword, then $w$ must be of type I; otherwise, $H(w)$ would contain one of the impermissible configurations of Lemma~\ref{lem:impermissible.heap.configs}.
\end{proof}

The next two lemmas are generalizations of Lemma 5.3 in~\cite{Green.R:P} and begin to describe the form that a non-cancellable element that is not of type I can take.  Recall from Section~\ref{subsec:heaps}, that $\r_k$ denotes the $k$th row of the canonical representation of the heap of a fully commutative element.

\begin{lem}\label{lem:weak.star.middle}
Let $w \in \FC(\C_{n})$ with $n \geq 4$ and suppose that $w$ is non-cancellable and not of type I.  If $s_{i} \in \r_{k+1}$ with $i \notin \{1,2,n,n+1\}$, then this entry is covered by entries labeled by $s_{i-1}$ and $s_{i+1}$.
\end{lem}

\begin{proof}
Note that our restrictions on $i$ and $n$ force $m(s_{i}, s_{i-1})=m(s_{i}, s_{i+1})=3$.  We proceed by induction.  

For the base case, assume that $k=1$.  Then the entry in $\r_2$ labeled by $s_i$ is covered by at least one of $s_{i-1}$ or $s_{i+1}$.   But since $w$ is not left weak star reducible, we must have both $s_{i-1}$ and $s_{i+1}$ occurring in $\r_{1}$.

For the inductive step, assume that the theorem is true for all $1 \leq k' \leq k-1$ for some $k$.  Suppose that $s_{i} \in \r_{k+1}$ with $i \notin \{1,2, n, n+1\}$.  Then at least one of $s_{i-1}$ or $s_{i+1}$ occur in $\r_{k}$.  We consider two cases: (1) $i \notin \{3,n-1\}$ and (2) $i  \in \{3,n-1\}$.

Case (1): Assume that in addition to $i \notin \{1,2,n,n+1\}$, $i \notin \{3,n-1\}$.  Observe that this forces $n \geq 6$.  Without loss of generality, assume that $s_{i-1} \in \r_{k}$.  (The case with $s_{i+1} \in \r_{k}$ is symmetric since the restrictions on $i$ imply that we may apply the induction hypothesis to either $i-1$ or $i+1$.)  By induction, the entry labeled by $s_{i-1}$ occurring in $\r_{k}$ is covered by an entry labeled by $s_{i-2}$ and an entry labeled by $s_{i}$.  This implies that the entry labeled by $s_{i}$ occurring in $\r_{k+1}$ must be covered by an entry labeled by $s_{i+1}$; otherwise, we produce one of the impermissible configurations of Lemma~\ref{lem:impermissible.heap.configs} corresponding to the subword $s_{i}s_{i-1}s_{i}$.  This yields our desired result.

Case (2):  For the second case, assume that $i \in \{3, n-1\}$.  Without loss of generality, assume that $i=3$; the case with $i=n-1$ can be handled by a symmetric argument.  Then $s_{3} \in \r_{k+1}$ and this entry is covered by either (a) an entry labeled by $s_{2}$, (b) an entry labeled by $s_{4}$, or (c) both.  If we are in situation (c), then we are done.  For sake of a contradiction, assume that exactly one of (a) or (b) occurs.

First, assume that (a) occurs, but (b) does not.  That is, $s_{2} \in \r_k$ and the entry labeled by $s_{3}$ that occurs in $\r_{k+1}$ is not covered by an entry labeled by $s_{4}$.  Since $k \geq 2$ and $w$ is fully commutative, it must be the case that $s_1\in\r_{k-1}$ while the entry labeled by $s_2$ occurring in $\r_k$ is not covered by an entry labeled by $s_3$.

For sake of a contradiction, assume that $k > 2$, so that $\r_{k-1}$ is not the top row of the canonical representation for $H(w)$.  Then we must have $s_{2} \in \r_{k-2}$.  Also, we cannot have $k=3$; otherwise, $w$ is left weak star reducible by $s_{2}$ with respect to $s_{1}$.  So, $k>3$, which implies that the entry labeled by $s_{2}$ occurring in $\r_{k-2}$ is covered.  This entry cannot be covered by $s_{1}$ since $w$ is fully commutative.  Therefore, we have $s_{3} \in \r_{k-3}$.  But by induction, this entry is covered by an entry labeled by $s_{2}$ and an entry labeled by $s_{4}$.  This produces one of the impermissible configurations of Lemma~\ref{lem:impermissible.heap.configs} corresponding to the subword $s_{2}s_{3}s_{2}$, which contradicts $w$ being fully commutative.  Thus, we must have $k=2$.  

Now, since $w$ is fully commutative and non-cancellable, we must conclude that 
\[
\includegraphics[scale=.85]{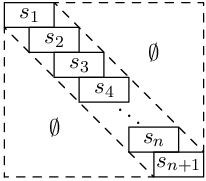}
\]
forms the top $n+1$ rows of the canonical representation of $H(w)$, where $\r_j=s_j$ for $1\leq j\leq n+1$.  Any other possibility either produces one of the impermissible configurations of Lemma~\ref{lem:impermissible.heap.configs} or violates $w$ not being right weak star reducible.  Since $w$ is not of type I, this cannot be all of $H(w)$.  The only possibility is that $s_n\in\r_{n+2}$.   Since $w$ is not right weak star reducible, this cannot be all of $H(w)$ either.  So, at least one of $s_{n-1}$ or $s_{n+1}$ occur in $\r_{n+3}$.  We cannot have $s_{n+1} \in \r_{n+3}$ because $w$ is fully commutative.  Thus, $s_{n-1} \in \r_{n+3}$ while $s_{n+1}$ is not.  Again, since $w$ is not right weak star reducible, we must have $s_{n-2} \in \r_{n+4}$ while $s_{n} \notin \r_{n+4}$.  Continuing with similar reasoning, we quickly see that
\[
\includegraphics[scale=.85]{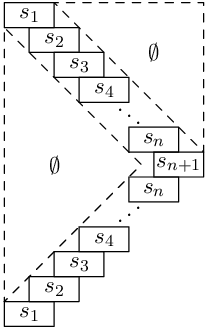}
\]
is a convex subheap of $H(w)$.  But then by Lemma~\ref{lem:zigzag}, $w$ is of type I, which is a contradiction.  Therefore, we cannot have possibility (a) occurring while (b) does not.

The only remaining possibility is that (b) occurs, but (a) does not.  That is, $s_{3} \in \r_{k+1}$ and $s_{4} \in \r_{k}$, while the entry labeled by $s_{3}$ occurring in $\r_{k+1}$ is not covered by an entry labeled by $s_{2}$.  Observe that the case $n=4$ is covered by an argument that is symmetric to the argument made above when we assumed that (a) occurs, but (b) does not, where we take $i=n-1$ instead of $i=3$.  So, assume that $n>4$.  Then by induction, entries labeled by $s_{3}$ and $s_{5}$ both cover $s_4\in\r_{k}$.  But then we produce one of the impermissible configurations of Lemma~\ref{lem:impermissible.heap.configs} corresponding to the subword $s_{3}s_{4}s_{3}$, which contradicts $w$ being fully commutative.  

We have exhausted all possibilities, and hence we have our desired result.
\end{proof}

\begin{lem}\label{lem:weak.star.end}
Let $w \in \FC(\C_{n})$ with $n\geq 3$ and suppose that $w$ is non-cancellable and not of type I.   If $s_{2}$ (respectively, $s_n$) occurs in $\r_{k+1}$ and is covered by an entry labeled by $s_{3}$ (respectively, $s_{n-1}$), then an entry labeled by $s_{1}$ (respectively, $s_{n+1}$) covers the entry labeled by $s_{2}$ (respectively, $s_n$) that occurs in $\r_{k+1}$.
\end{lem}

\begin{proof}
Assume that $s_2$ occurs in $\r_{k+1}$ and is covered by an entry labeled by $s_3$.  The case involving $s_n$ being covered by an entry labeled by $s_{n-1}$ follows by a symmetric argument.  The case with $n=3$ proves to be more difficult than when $n>3$.  We handle this more difficult case first.

Assume that $n=3$.  If $k=1$, then $s_{2} \in \r_{2}$ and $s_{3} \in \r_{1}$.    This implies that $s_{1}$ must occur in $\r_{1}$; otherwise, $w$ is left weak star reducible by $s_{3}$ with respect to $s_{2}$.  So, assume that $k \geq 2$.  For sake of a contradiction, assume that the entry labeled by $s_{2}$ occurring in $\r_{k+1}$ is not covered by an entry labeled by $s_{1}$.  This forces $s_{3} \in \r_{k}$.  Then at least one of $s_{2}$ or $s_{4}$ must cover the entry labeled by $s_{3}$ occurring in $\r_{k-1}$.  Since $s_{1}$ does not cover the occurrence of $s_{2} \in \r_{k+1}$ and $w$ is fully commutative, it must be the case that $s_{4} \in \r_{k-1}$, while the entries labeled by $s_{3}$ and $s_{2}$ occurring in $\r_{k}$ and $\r_{k+1}$, respectively, are not covered by entries labeled by $s_{2}$ and $s_{1}$, respectively.  Then
\[
\includegraphics[scale=.85]{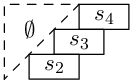}
\]
is a convex subheap of $H(w)$.

First, assume that $k=2$, so that $s_{2} \in \r_{3}$, $s_{3} \in \r_{2}$, $s_{4} \in \r_{1}$, and neither $s_{1}$ nor $s_{2}$ occur in $\r_{1}$ or $\r_{2}$.  Then the subheap immediately above is the northwest corner of $H(w)$, where the entry labeled by $s_{4}$ occurs in the top row.  Since $w$ is not of type I, this cannot be all of $H(w)$.  Furthermore, since $w$ is fully commutative and not right weak star reducible, it must be the case that 
\[
\includegraphics[scale=.85]{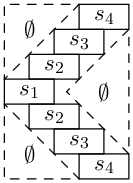}
\]
is the top of $H(w)$.  But then by Lemma~\ref{lem:zigzag}, $w$ is of type I, which is a contradiction.

Next, assume that $k>2$.  In this case, the entry in $\r_{k-1}$ labeled by $s_{4}$ must be covered by an entry labeled by $s_{3}$.  This implies that $s_{4}$ cannot occur in $\r_{k+1}$ since $w$ is fully commutative.  However, since $w$ is fully commutative and not left weak star reducible, it must be the case that 
\[
\includegraphics[scale=.85]{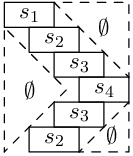}
\]
is a convex subheap of $H(w)$.  Since $w$ is not of type I, this cannot be all of $w$.  Since $w$ is fully commutative, the only two possibilities are that the entry labeled by $s_{1}$ in $\r_{k-4}$ is covered by an entry labeled by $s_{2}$ or that the entry labeled by $s_{2}$ in $\r_{k+1}$ covers an entry labeled by $s_{1}$.  In either case, $w$ is of type I by Lemma~\ref{lem:zigzag}, which yields a contradiction.

We have exhausted all possibilities.  Therefore, it must be the case that the entry labeled by $s_{2} \in \r_{k+1}$ is covered by an entry labeled by $s_{1}$, as desired.  This completes the case when $n=3$.

Now, assume that $n\geq 4$.  Note that since $n\geq 4$, $m(s_{2},s_{3})=3$.  If $k=1$, then $s_{2} \in \r_{2}$ and $s_{3} \in \r_{1}$.    This implies that $s_{1}$ must occur in $\r_{1}$, otherwise, $w$ is left weak star reducible by $s_{3}$ with respect to $s_{2}$.  Assume that $k \geq 2$, so that $s_{2} \in \r_{k+1}$ is covered by an entry labeled by $s_{3}$.  Then by Lemma~\ref{lem:weak.star.middle}, entries labeled by $s_{2}$ and $s_{4}$ cover the entry labeled by $s_{3}$ occurring in $\r_{k}$.  Since $w$ is fully commutative, we must have the entry labeled by $s_{2}$ occurring in $\r_{k+1}$ covered by an entry labeled by $s_{1}$, as desired; otherwise, we produce one of the impermissible configurations of Lemma~\ref{lem:impermissible.heap.configs} and violate $w$ being fully commutative.  This completes the case when $n\geq 4$.
\end{proof}

\begin{rem}
Lemmas~\ref{lem:weak.star.middle} and~\ref{lem:weak.star.end} have ``upside-down'' versions, where we replace $k+1$ with $k-1$ and we swap the phrases ``is covered by'' and ``covers.'' 
\end{rem}

The next four lemmas are all of a similar flavor.  In each case, we require $w \in \FC(\C_n)$ to be non-cancellable and we have one lemma for each of the following rank situations: (1) $n=2$, (2) $n=3$, (3) $n=4$, and (4) $n>4$.  We have divided the lemmas up into these four situations because the respective proofs are different and in some cases we require slight modifications to the hypotheses.  These lemmas will be used to construct elements of type II in the proof of the classification of the non-cancellable elements of type $\C_n$.

\begin{lem}\label{lem:sandwich.stack.n=2}
Let $w \in \FC(\C_{2})$ and suppose that $w$ is non-cancellable and not of type I.  If $\r_{k+1}=\x_{\O}$ (respectively, $\x_{\E}$), then $\r_{k}=\x_{\E}$ (respectively, $\x_{\O}$).
\end{lem}

\begin{proof}
Note that when $n=2$, we have $\x_{\O}=s_{1}s_{3}$ and $\x_{\E}=s_{2}$.  If $\r_{k+1}=s_{1}s_{3}$, then it is clear that $\r_{k}=s_{2}$.  Now, assume that $\r_{k+1}=s_{2}$.  Then at least one of $s_{1}$ or $s_{3}$ occurs in $\r_{k}$.  For sake of a contradiction, assume that only one of these occurs in $\r_{k}$, and without loss of generality, assume that $s_{1} \in \r_{k}$ while $s_{3} \notin \r_{k}$; the remaining case is handled by a symmetric argument.  We consider two cases: (1) $k=1$ and (2) $k \geq 2$.

Case (1):  First, assume that $k=1$, so that $s_{1} \in \r_{1}$ and $s_{2} \in \r_{2}$, while $s_{3} \notin \r_{1}$.  This cannot be all of $H(w)$ since $w$ is not of type I.  Since $w$ is fully commutative and non-cancellable, the only possibility is that 
\[
\includegraphics[scale=.85]{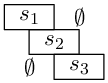}
\]
forms the top three rows of the canonical representation of $H(w)$.  Again, since $w$ is not of type I, there must be more to $H(w)$.  The top five rows of $H(w)$ must equal
\[
\includegraphics[scale=.85]{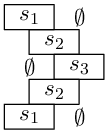}
\]
since $w$ is fully commutative and non-cancellable.  But then according to Lemma~\ref{lem:zigzag}, $w$ is of type I, which is a contradiction.

Case (2):  For the second case, assume that $k\geq 2$.  Then we must have $s_{2} \in \r_{k-1}$.  Since $w$ is not left weak star reducible, we cannot have $k=2$; otherwise, $w$ is left weak star reducible by $s_{2}$ with respect to $s_{1}$.  Thus, $k>2$, and hence at least one of $s_{1}$ or $s_{3}$ occurs in $\r_{k-2}$.  Since $w$ is fully commutative, $s_{1} \notin \r_{k-2}$, and so, $s_{3} \in \r_{k-2}$.  This implies that
\[
\includegraphics[scale=.85]{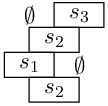}
\]
is a convex subheap of $H(w)$.  This cannot be all of $H(w)$ since $w$ is not of type I.  The only possibilities are that $s_{2} \in \r_{k-3}$ or $s_{3} \in \r_{k+2}$ (both possibilities could occur simultaneously).  In either case, $w$ must be of type I by Lemma~\ref{lem:zigzag}, which is a contradiction.
\end{proof}

Unlike the previous lemma, the next lemma does not place any requirements on whether $w$ is or is not of type I.

\begin{lem}\label{lem:sandwich.stack.n=3}
Let $w \in \FC(\C_{3})$ and suppose that $w$ is non-cancellable.  If $\r_{k+1}=\x_{\O}$ (respectively, $\x_{\E}$), then $\r_{k}=\x_{\E}$ (respectively, $\x_{\O}$).
\end{lem}

\begin{proof}
Note that when $n=3$, we have $\x_{\O}=s_{1}s_{3}$ and $\x_{\E}=s_{2}s_{4}$.  Assume that $\r_{k+1}=\x_{\O}$; the proof of the other case is similar.  Then we must have $s_{2} \in \r_{k}$ since this is the only generator available to cover $s_{1} \in \r_{k+1}$.  Since $\r_{k+1}=\x_{\O}$, $n(w)>1$, which implies that $w$ is not of type I by Proposition~\ref{prop:zigzags}.  Then by Lemma~\ref{lem:weak.star.end}, the entry labeled by $s_{3} \in \r_{k+1}$ must be covered by an entry labeled by $s_{4}$.  If $k=1$, then $s_{4} \in \r_{k}$, as desired.  If $k>1$, then at least one of $s_{1}$ or $s_{3}$ occurs in $\r_{k-1}$.  For sake of a contradiction, assume that $s_{1} \in \r_{k-1}$, but $s_{3} \notin \r_{k-1}$.  If $k=2$, then $w$ would be left weak star reducible by $s_{1}$ with respect to $s_{2}$.  So, we must have $k>2$, in which case, $s_{2} \in \r_{k-2}$.  But then we produce one of the impermissible configurations of Lemma~\ref{lem:impermissible.heap.configs} corresponding to the subword $s_{2}s_{1}s_{2}s_{1}$, which contradicts $w$ being fully commutative.  So, it must be the case that  $s_{3} \in \r_{k-1}$.  This implies that the entry labeled by $s_{4}$ that covers $s_{3} \in \r_{k+1}$ must occur in $\r_{k}$.  So, $\r_{k}=\x_{\E}$, as desired.
\end{proof}

Notice that the next lemma requires $w$ to have full support, which the previous two lemmas did not.  (Recall that $w$ has full support if $\supp(w)=S$; see Section~\ref{subsec:Coxeter_groups} for the definition.)

\begin{lem}\label{lem:sandwich.stack.n=4}
Let $w \in \FC(\C_{4})$ and suppose that $w$ is non-cancellable and has full support.  If $\r_{k+1}=\x_{\O}$ (respectively, $\x_{\E}$), then $\r_{k}=\x_{\E}$ (respectively, $\x_{\O}$).
\end{lem}

\begin{proof}
Note that when $n=4$, we have $\x_{\O}=s_{1}s_{3}s_{5}$ and $\x_{\E}=s_{2}s_{4}$.   We consider two cases: (1) $\r_{k+1}=\x_{\O}$ and (2) $\r_{k+1}=\x_{\E}$.

Case (1):  First, assume that $\r_{k+1}=\x_{\O}$.  By Lemma~\ref{lem:weak.star.middle}, the entry labeled by $s_{3} \in \r_{k+1}$ is covered by entries labeled by $s_{2}$ and $s_{4}$, where at least one of these occurs in $\r_{k}$.  Since $s_{2}$ (respectively, $s_{4}$) is the only generator that may cover an entry labeled by $s_{1}$ (respectively, $s_{5}$), we must have both $s_{2}$ and $s_{4}$ occurring in $\r_{k}$, as desired.

Case (2):  For the more difficult case, assume that $\r_{k+1}=\x_{\E}$.  First, we argue that an entry labeled by $s_{3}$ covers the occurrences of $s_{2}$ and $s_{4}$ in $\r_{k+1}$.  For sake of a contradiction, assume otherwise.  Then we must have $s_{1}$ and $s_{5}$ both occurring in $\r_{k}$ since these are the only entries available to cover the entries occurring in $\r_{k+1}$.  We consider to subcases: (a) $k=1$ and (b) $k>2$.  

(a) Assume that $k=1$, so that $\r_{1}=s_{1}s_{5}$ and $\r_{2}=s_{2}s_{4}$.  We cannot have $s_{1}$ (respectively, $s_{5}$) occurring in $\r_{3}$; otherwise $w$ would be left weak star reducible by $s_{1}$ (respectively, $s_{5}$) with respect to $s_{2}$ (respectively, $s_{4}$).  Since $w$ has full support, we must have $s_{3} \in \r_{3}$.  By the upside-down version of  Lemma~\ref{lem:weak.star.middle}, the entry labeled by $s_{3} \in \r_{3}$ must cover entries labeled by $s_{2}$ and $s_{4}$.  But this produces impermissible configurations of Lemma~\ref{lem:impermissible.heap.configs} corresponding to the subwords $s_{2}s_{3}s_{2}$ and $s_{4}s_{3}s_{4}$, which contradicts $w$ being fully commutative.  

(b) Next, assume that $k\geq 2$.  Then we must have $s_{2}$ and $s_{4}$ occurring in $\r_{k-1}$.  Since $w$ is not left weak star reducible, we must have $k > 2$; otherwise, $w$ is left weak star reducible by $s_{2}$ (respectively, $s_{4}$) with respect to $s_{1}$ (respectively, $s_{5}$).  The entry labeled by $s_{2}$ (respectively, $s_{4}$) occurring in $\r_{k-1}$ cannot be covered by $s_{1}$ (respectively, $s_{5}$); otherwise, we produce one of the impermissible configurations of Lemma~\ref{lem:impermissible.heap.configs}.  So, we must have $s_{3} \in \r_{k-2}$.  If $k=3$, then $w$ would be left weak star reducible by $s_{3}$ with respect to either $s_{2}$ or $s_{4}$.  Thus, $k \geq 4$.  By Lemma~\ref{lem:weak.star.middle}, the entry labeled by $s_{3} \in \r_{k-2}$ is covered by entries labeled by $s_{2}$ and $s_{4}$.  But then we again produce impermissible configurations of Lemma~\ref{lem:impermissible.heap.configs} corresponding to the subwords $s_{2}s_{3}s_{2}$ and $s_{4}s_{3}s_{4}$, which is a contradiction.

We have shown that if $\r_{k+1}=\x_{\E}$, then the entries labeled by $s_{2}$ and $s_{4}$ occurring in $\r_{k+1}$ must be covered by an entry labeled by $s_{3}$.  By Lemma~\ref{lem:weak.star.end}, an entry labeled by $s_{1}$ (respectively, $s_{5}$) covers the entry labeled by $s_{2}$ (respectively, $s_{4}$) occurring in $\r_{k+1}$.   If $\r_{k} \neq s_{1}s_{3}s_{5}$, we quickly contradict Lemma~\ref{lem:impermissible.heap.configs} or Lemma~\ref{lem:weak.star.middle}.  Therefore, we must have $\r_{k}=\x_{\O}$, as desired.
\end{proof}

We have reached the last of our preparatory lemmas.

\begin{lem}\label{lem:sandwich.stack.n>4}
Let $w \in \FC(\C_{n})$ with $n > 4$ and suppose that $w$ is non-cancellable.  If $\r_{k+1}=\x_{\O}$ (respectively, $\x_{\E}$), then $\r_{k}=\x_{\E}$ (respectively, $\x_{\O}$).
\end{lem}

\begin{proof}
If $k=1$, then the result follows by Lemmas~\ref{lem:weak.star.middle} and~\ref{lem:weak.star.end}.  If $k>1$, the result follows by making repeated applications of Lemmas~\ref{lem:weak.star.middle} and~\ref{lem:weak.star.end} while avoiding the impermissible configurations of Lemma~\ref{lem:impermissible.heap.configs}.
\end{proof}

\begin{rem}
Lemmas~\ref{lem:sandwich.stack.n=2}--\ref{lem:sandwich.stack.n>4} all have ``upside-down'' versions since all of the arguments reverse nicely.  That is, if $w$ is non-cancellable (and not of type I when $n=2$; and has full support when $n=4$), then $\r_{k}=\x_{\O}$ (respectively, $\x_{\E}$) implies $\r_{k+1}=\x_{\E}$ (respectively, $\x_{\O}$).
\end{rem}

\end{subsection}

%%%%%%%%%%%%%%%%%%%%%

\begin{subsection}{Proof of classification}

We are now ready to prove the classification of the type $\C$ non-cancellable elements.

%\vspace{1em}

%I screwed this up when I submitted to the arXiv!
\begin{proof}[Proof of Theorem~\ref{thm:affineCwsrm}]
%\noindent \textbf{Proof of Theorem~\ref{thm:affineCwsrm}.}  
It is easily seen that every element on our list is non-cancellable.  For sake of a contradiction, assume that there exists $w \in \FC(\C_{n})$ such that $w$ is non-cancellable, but not on our list.  If there exists $s \notin \supp(w)$, then $w$ is equal to $uv$ (reduced), where $u$ is of type $B$, $v$ is of type $B'$, and $\supp(u) \cap \supp(v) = \emptyset$.  Since $w$ is non-cancellable, both $u$ and $v$ are non-cancellable.  But then $w$ must be one of the elements from (i), which contradicts our assumption that $w$ is not on our list.  So, if $w$ is not on our list, $w$ must have full support.  In particular, this implies that $w$ is not a product of commuting generators.  According to Proposition~\ref{prop:zigzags}, the only non-cancellable elements with $n(w)=1$ are already listed in (ii).  Hence $n(w)>1$ (i.e., $w$ is not of type I).  

Now, consider the canonical representation of $H(w)$ and suppose that it has $m$ rows.  Since $w$ is not a product of commuting generators, $m\geq 2$.  Our immediate goal is to show that $\r_m$ is equal to either $\x_\O$ or $\x_\E$.  Then we will be able to make use of Lemmas~\ref{lem:sandwich.stack.n=2}--\ref{lem:sandwich.stack.n>4} to conclude that $w$ is of type II.  We consider three main cases: (1) $n=2$, (2) $n=3$, and (3) $n \geq 4$.

Case (1):  Assume that $n=2$.  In this case, $\x_{\O}=s_{1}s_{3}$ and $\x_{\E}=s_{2}$.  If $s_{2} \in \r_{m}$, then $\r_{m}=\x_{\E}$.  Assume that $s_{2} \notin \r_{m}$.  Then at least one of $s_{1}$ or $s_{3}$ occurs in $\r_{m}$.  Then we must have $s_{2} \in \r_{m-1}$.  In fact, $\r_{m-1}=\x_{\E}$.  By the upside-down version of Lemma~\ref{lem:sandwich.stack.n=2}, $\r_{m}=\x_{\O}$.

Case (2):  For the second case, assume that $n=3$.  In this case, $\x_{\O}=s_{1}s_{3}$ and $\x_{\E}=s_{2}s_{4}$.  For sake of a contradiction, assume that $s_{i} \in \r_{m}$ but $s_{i'} \notin \r_{m}$, where $|i-i'|=2$.  We consider the subcases: (a) $i=1$ and (b) $i=3$.  The cases $i=2$ and $i=4$ are similar.

(a)  Suppose that $i=1$, so that $s_{1} \in \r_{m}$ while $s_{3} \notin \r_{m}$.  Then we must have $s_{2} \in \r_{m-1}$.  Since $w$ has full support and $s_{3}$ does not occur in $\r_{m-1}$ or $\r_{m}$, we must have $m \geq 3$.  Then the entry labeled by $s_{2}$ occurring in $\r_{m-1}$ cannot be covered by an entry labeled by $s_{1}$; otherwise, $w$ is right weak star reducible by $s_{1}$ with respect to $s_{2}$.  Thus, $s_{3} \in \r_{m-2}$.  But according to Lemma~\ref{lem:weak.star.end}, the entry labeled by $s_{2}$ occurring in $\r_{m-1}$ must be covered by an entry labeled by $s_{1}$, which is a contradiction.  

(b)  Next, suppose that $i=3$, so that $s_{3} \in \r_{m}$ while $s_{1} \notin \r_{m}$.  Then at least one of $s_{2}$ or $s_{4}$ occurs in $\r_{m-1}$.  If $s_{2} \in \r_{m-1}$, then $w$ would be right weak star reducible by $s_{3}$ with respect to $s_{2}$.  So, it must be the case that $s_{4} \in \r_{m-1}$, while the entry labeled by $s_{3}$ is not covered by an entry labeled by $s_{2}$.  Since $w$ has full support, we must have $m\geq 3$; otherwise, $s_{1}, s_{2} \notin \supp(w)$.  This implies that $s_{3} \in \r_{m-2}$.  But then $w$ is right weak star reducible by $s_{3}$ with respect to $s_{4}$, which is again a contradiction.

Case (3):  Lastly, assume that $n\geq 4$.  For sake of a contradiction, assume that $s_{i} \in \r_{m}$ but $s_{i'} \notin \r_{m}$, where $|i-i'|=2$.  Without loss of generality, assume that $1\leq i \leq n-1$ and $i'=i+2$, so that $s_{i+2} \notin \r_{m}$; the remaining cases are similar.  We consider three possibilities: (a) $i=1$, (b) $i=2$, and (c) $3\leq i \leq n-1$.

(a) If $i=1$, then this case is identical to (a) in Case (2), where we contradict Lemma~\ref{lem:weak.star.end}.

(b) Next, assume that $i=2$, so that $s_{2} \in \r_{m}$ while $s_{4} \notin \r_{m}$.  Then the entry labeled by $s_{2} \in \r_{m}$ cannot be covered by an entry labeled by $s_{3}$; otherwise, $w$ would be right weak star reducible by $s_{2}$ with respect to $s_{3}$.  This implies that we must have $s_{1} \in \r_{m-1}$.  Since $w$ has full support and $s_{3}$ does not occur in $\r_{m-1}$ or $\r_{m}$, we must have $m \geq 3$.  Then $s_{2} \in \r_{m-2}$.  But then $w$ is right weak star reducible by $s_{2}$ with respect to $s_{1}$, which is a contradiction.

(c) Finally, assume that $3\leq i \leq n-1$, so that $s_{i} \in \r_{m}$ while $s_{i+2} \notin \r_{m}$.  By Lemma~\ref{lem:weak.star.middle}, the entry labeled by $s_{i} \in \r_{m}$ must be covered by entries labeled by $s_{i-1}$ and $s_{i+1}$.  However, this implies that $w$ is right weak star reducible by $s_{i}$ with respect to $s_{i+1}$, which is a contradiction.

In any case, we have shown that $\r_m$ is equal to either $\x_\O$ or $\x_\E$.  By making repeated applications of Lemma~\ref{lem:sandwich.stack.n=2} if $n=2$,  Lemma~\ref{lem:sandwich.stack.n=3} if $n=3$, Lemma~\ref{lem:sandwich.stack.n=4} if $n=4$, or Lemma~\ref{lem:sandwich.stack.n>4} if $n>4$, $w$ must be equal to an alternating product of $\x_{\O}$ and $\x_{\E}$.  This implies that $w$ is of type II, which contradicts our assumption that $w$ is not on our list. %\hfill $\qed$
\end{proof}

\end{subsection}

\end{section}

%%%%%%%%%%%%%%%%%%%%%

\begin{section}{Closing remarks and further research}\label{sec:closing}

The (type $A$) Temperley--Lieb algebra $\TL(A)$, invented by H.N.V.~Temperley and E.H.~Lieb in 1971~\cite{Temperley.H;Lieb.E:A}, is a finite dimensional associative algebra which arose in the context of statistical mechanics.  A diagram algebra is an associative algebra with a basis given by certain diagrams, in which the multiplication rule in the algebra is given by applying local combinatorial rules to the diagrams.  R.~Penrose and L.H.~Kauffman showed that $\TL(A)$ can be faithfully represented by a diagram algebra.

In 1987, V.F.R.~Jones showed that the (type $A$) Temperley--Lieb algebra occurs naturally as a quotient of the type $A$ Hecke algebra, $\H(A)$~\cite{Jones.V:B}.  If $(W,S)$ is Coxeter system of type $\Gamma$, the associated Hecke algebra $\H(\Gamma)$ is an algebra with a basis given by $\{T_w: w \in W\}$ and relations that deform the relations of $W$ by a parameter $q$. The realization of the Temperley--Lieb algebra as a Hecke algebra quotient was generalized by J.J.~Graham in~\cite{Graham.J:A} to the case of an arbitrary Coxeter system, which we denote by $\TL(\Gamma)$.  

One motivation behind studying these generalized Temperley--Lieb algebras is that they provide a gateway to understanding the Kazhdan--Lusztig theory of the associated Hecke algebra.  Loosely speaking, $\TL(\Gamma)$ retains some of the relevant structure of $\H(\Gamma)$, yet is small enough that computation of the leading coefficients of the notoriously difficult to compute Kazhdan--Lusztig polynomials is often much simpler.

Since Coxeter groups of type $\C$ have an infinite number of fully commutative elements, $\TL(\C)$ is infinite dimensional.  With the exception of type $\widetilde{A}$, all other generalized Temperley--Lieb algebras with known diagrammatic representations are finite dimensional.  In the finite dimensional case, counting arguments are employed to prove faithfulness, but these techniques are not available in the type $\C$ case.  The classification of the non-cancellable elements in Theorem~\ref{thm:affineCwsrm} provides the foundation for inductive arguments used to prove the faithfulness of the diagram algebra introduced by the author in~\cite{Ernst.D:A}.  This diagram algebra is the first faithful representation of an infinite dimensional non-simply laced generalized Temperley--Lieb algebra (in the sense of Graham).  Chapters 6--10 of the author's PhD thesis are concerned with establishing this representation and will be the focus of subsequent papers.  In a future paper, we plan to construct a Jones-type trace on $\H(\C)$ using the diagrammatic representation of $\TL(\C)$, allowing us to non-recursively compute leading coefficients of Kazhdan--Lusztig polynomials indexed by pairs of fully commutative elements.

\end{section}

%%%%%%%%%%%%%%%%%%%%%

\section*{Acknowledgements}

I would like to thank R.M.~Green for many useful conversations during the preparation of this article.  I am also grateful to the referee for his or her careful reading of the paper and constructive suggestions for improvements.

%%%%%%%%%%%%%%%%%%%%%

\bibliographystyle{plain}

\bibliography{non-cancellable}

\begin{thebibliography}{10}

\bibitem{Billey.S;Jones.B:A}
S.C. Billey and B.C. Jones.
\newblock Embedded factor patterns for {D}eodhar elements in
  {K}azhdan--{L}usztig theory.
\newblock {\em Ann. Comb.}, 11(3--4):285--333, 2007.

\bibitem{Billey.S;Warrington.G:A}
S.C. Billey and G.S. Warrington.
\newblock Kazhdan--{L}usztig polynomials for 321-hexagon-avoiding permutations.
\newblock {\em J. Algebraic Combin.}, 13:111--136, 2001.

\bibitem{Cartier.P;Foata.D:A}
P.~Cartier and D.~Foata.
\newblock Probl\`emes combinatoires de commutation et r\'earrangements.
\newblock {\em Lecture Notes in Mathematics, Springer-Verlag, New York/Berlin},
  85, 1969.

\bibitem{Ernst.D:A}
D.C. Ernst.
\newblock {\em A diagrammatic representation of an affine ${C}$
  {T}emperley--{L}ieb algebra}.
\newblock PhD thesis, University of Colorado at Boulder, 2008.

\bibitem{Fan.C:A}
C.K. Fan.
\newblock Structure of a {H}ecke algebra quotient.
\newblock {\em J. Amer. Math. Soc.}, 10:139--167, 1997.

\bibitem{Fan.C:B}
C.K. Fan.
\newblock Schubert varieties and short braidedness.
\newblock {\em Transform. Groups}, 3(1):51--56, 1998.

\bibitem{Fan.C;Green.R:A}
C.K. Fan and R.M. Green.
\newblock On the affine {T}emperley--{L}ieb algebras.
\newblock {\em Jour. L.M.S.}, 60:366--380, 1999.

\bibitem{Geck.M;Pfeiffer.G:A}
M.~Geck and G.~Pfeiffer.
\newblock {\em Characters of finite {C}oxeter groups and {I}wahori--{H}ecke
  algebras}.
\newblock Oxford University Press, 2000.

\bibitem{Graham.J:A}
J.J. Graham.
\newblock {\em Modular representations of {H}ecke algebras and related
  algebras}.
\newblock PhD thesis, University of Sydney, 1995.

\bibitem{Green.R:P}
R.M. Green.
\newblock Star reducible {C}oxeter groups.
\newblock {\em Glasgow Math. J.}, 48:583--609, 2006.

\bibitem{Green.R:K}
R.M. Green.
\newblock Generalized {J}ones traces and {K}azhdan--{L}usztig bases.
\newblock {\em J. Pure Appl. Alg.}, 211:744--772, 2007.

\bibitem{Green.R;Losonczy.J:E}
R.M. Green and J.~Losonczy.
\newblock A projection property for {K}azhdan--{L}usztig bases.
\newblock {\em Int Math Res Notices}, 1:23--34, 2000.

\bibitem{Green.R;Losonczy.J:B}
R.M. Green and J.~Losonczy.
\newblock Fully commutative {K}azhdan--{L}usztig cells.
\newblock {\em Ann. Inst. Fourier Grenoble}, 51:1025--1045, 2001.

\bibitem{Humphreys.J:A}
J.E. Humphreys.
\newblock {\em Reflection Groups and Coxeter Groups}.
\newblock Cambridge University Press, 1990.

\bibitem{Jones.V:B}
V.F.R. Jones.
\newblock Hecke algebra representations of braid groups and link polynomials.
\newblock {\em Ann. of Math. 2}, 126:335--388, 1987.

\bibitem{Kazhdan.D;Lusztig.G:A}
D.~Kazhdan and G.~Lusztig.
\newblock Representations of {C}oxeter groups and {H}ecke algebras.
\newblock {\em Invent. Math.}, 53:165--184, 1979.

\bibitem{Lusztig.G:A}
G.~Lusztig.
\newblock Cells in affine {W}eyl groups, {I}.
\newblock In {\em Algebraic groups and related topics}, pages 255--287. Adv.
  Studies Pure Math 6, North-Holland and Kinokuniya, Tokyo and Amsterdam, 1985.

\bibitem{Shi.J:B}
J.Y. Shi.
\newblock Fully commutative elements and {K}azhdan--{L}usztig cells in the
  finite and affine {C}oxeter groups.
\newblock {\em Proc. Amer. Math. Soc.}, 131:3371--3378, 2003.

\bibitem{Shi.J:C}
J.Y. Shi.
\newblock Fully commutative elements in the {W}eyl and affine {W}eyl groups.
\newblock {\em J. Algebra}, 284:13--36, 2005.

\bibitem{Stembridge.J:B}
J.R. Stembridge.
\newblock On the fully commutative elements of {C}oxeter groups.
\newblock {\em J. Algebraic Combin.}, 5:353--385, 1996.

\bibitem{Stembridge.J:A}
J.R. Stembridge.
\newblock The enumeration of fully commutative elements of {C}oxeter groups.
\newblock {\em J. Algebraic Combin.}, 7(3):291--320, 1998.

\bibitem{Temperley.H;Lieb.E:A}
H.N.V. Temperley and E.H. Lieb.
\newblock Relations between percolation and colouring problems and other graph
  theoretical problems associated with regular planar lattices: some exact
  results for the percolation problem.
\newblock {\em Proc. Roy. Soc. London Ser. A}, 322:251--280, 1971.

\bibitem{Viennot.G:A}
G.X. Viennot.
\newblock Heaps of pieces, {I}: basic definitions and combinatorial lemmas.
\newblock In G.~Labelle and P.~Leroux, editors, {\em Combinatoire
  \'E\-nu\-m\'e\-ra\-tive}, pages 321--350. Springer-Verlag, 1986.

\end{thebibliography}

\end{document}